\documentclass[9pt]{amsart}
\usepackage{graphicx}
\usepackage{amssymb}
\usepackage{amsmath}
\usepackage{amsthm}
\usepackage{amscd}
\usepackage[all,2cell,dvips]{xy}
\UseAllTwocells \SilentMatrices
\newtheorem{thm}{Theorem}[section]

\newtheorem{cor}[thm]{Corollary}
\newtheorem{lem}[thm]{Lemma}

\newtheorem{prop}[thm]{Proposition}
\theoremstyle{definition}

\theoremstyle{remark}
\newtheorem{rem}[thm]{\bf Remark}
\numberwithin{equation}{section}

\newcommand{\A}{\Lambda}

\begin{document}

\sf

\title[Calabi-Yau objects in triangulated categories]{Calabi-Yau objects in triangulated categories}
\author[C. Cibils, P. Zhang] {Claude Cibils$^a$ and Pu Zhang$^{b, *}$}
\thanks{The second named author is supported by the CNRS of France, the NSF of
China and of Shanghai City (Grant No. 10301033 and ZR0614049).}
\thanks{$^*$ The corresponding author}
\maketitle

\begin{center}
$^A$D\'epartement de Math\'ematiques, \ \ Universit\'e de Montpellier 2\\
F-34095, Montpellier Cedex 5, France\ \ \
Claude.Cibils$\symbol{64}$math.univ-montp2.fr\\
$^B$Department of Mathematics, \ \ Shanghai Jiao Tong University\\
Shanghai 200240, P. R. China\ \ \ \ pzhang$\symbol{64}$sjtu.edu.cn
\end{center}
\begin{abstract} We introduce the Calabi-Yau (CY) objects in a
Hom-finite Krull-Schmidt triangulated $k$-category, and notice
that the structure of the minimal, consequently all the CY
objects, can be described. The relation between indecomposable CY
objects and Auslander-Reiten triangles is provided. Finally we
classify all the CY modules of self-injective Nakayama algebras,
determining this way the self-injective Nakayama algebras
admitting indecomposable CY modules. In particular, this result
recovers the algebras whose stable categories are Calabi-Yau,
which have been obtained in [BS].
\end{abstract}

\vskip10pt

\section {\bf Introduction}

Calabi-Yau (CY) categories have been introduced by Kontsevich
[Ko]. They provide a new insight and a wide framework for topics
as in mathematical physics ([Co]), non-commutative geometry
([B], [Gin1], [Gin2]), and representation theory of Artin algebras
([BS], [ES], [IR], [Ke], [KR1], [KR2]).

Triangulated categories with Serre dualities ([BK], [RV]) and CY
categories have important global naturality. On the other hand,
even in non CY categories, inspired by [Ko], one can introduce CY
objects. It turns out that they arise naturally in non CY
categories and enjoy ``local naturality'' and interesting
properties (Prop. 4.4, Theorems 3.2, 4.2, 5.5 and 6.1).

\vskip10pt

The first aim of this paper is to study the properties of such
objects in a Hom-finite Krull-Schmidt triangulated $k$-category
with Serre functor $F$. We give the relation between
indecomposable CY objects and the Auslander-Reiten triangles ($\S
3$), and describe all the $d$-th CY objects via the minimal ones,
which are exactly the direct sum of all the objects in finite
$\langle [-d]\circ F\rangle$-orbits of
$\operatorname{Ind}(\mathcal A)$ ($\S 4$). We classify all the
$d$-th CY modules of self-injective Nakayama algebras for any
integer $d$ ($\S 5$). Finally, we determine all the self-injective
Nakayama algebras which admit indecomposable CY modules. In
particular, this recovers the algebras whose stable categories are
Calabi-Yau ($\S 6$), included in the work of Bialkowski and
Skowro\'nski [BS]. Note that the CY modules are invariant under
stable equivalences between self-injective algebras, with a very
few exceptions (Prop.3.1). Consequently our results on
self-injective Nakayama algebras extend to the one on the
wreath-like algebras ([GR]), which contains the Brauer tree
algebras ([J]).

This also raises an immediate question. Let $\mathcal A$ be a
Hom-finite Krull-Schmidt triangulated $k$-category  with a Serre
functor. If all objects are $d$-th CY with the same $d$, whether
$\mathcal A$ is a Calabi-Yau category?

\vskip10pt

\section {\bf Backgrounds and Preliminaries}

\subsection{} Let $k$ be a field and $\mathcal{A}$ a Hom-finite
$k$-category. Recall from Bondal and Kapranov [BK]
that a $k$-linear functor $F \colon
\mathcal{A} \to \mathcal{A}$ \ is {\em a right Serre functor} if
there exist $k$-isomorphisms $$\eta_{A, B}: \ \
\operatorname{Hom}_{\mathcal{A}}(A, B) \longrightarrow  D
\operatorname{Hom}_{\mathcal{A}}(B, FA), \ \ \forall \ \ A, \ B\in
\mathcal{A},$$ which are natural both in $A$ and $B$, where $D =
\operatorname{Hom}_{k}(-, k)$. Such an $F$ is unique up to a natural isomorphism, and fully-faithful;
if it is an equivalence, then a quasi-inverse $F^{-1}$ is a left Serre
functor; in this case we call $F$ a Serre functor. Note that
$\mathcal{A}$ has a Serre functor if and only if it has both right
and left Serre functor. See Reiten and Van den Bergh [RV].

\vskip10pt

For triangulated categories we refer to [Har], [V], and [N]. Let
$\mathcal{A}$ be a Hom-finite triangulated $k$-category. Following
Happel [Hap1], {\em an Auslander-Reiten triangle} $X\stackrel{f}
{\longrightarrow} Y \stackrel{g} {\longrightarrow} Z \stackrel{h}
{\longrightarrow} X[1]$ of $\mathcal{A}$ is a distinguished
triangle satisfying:

(AR1) \ $X$ and $Z$ are indecomposable;

(AR2) \ $h\ne 0$;

(AR3) \ If $t: Z'\longrightarrow Z$ is not a retraction, then
there exists $t': Z'\longrightarrow Y$ such that $t = gt'$.

\vskip10pt

Note that (AR3) is equivalent to

(AR4) \ If $Z'$ is indecomposable and $t: Z'\longrightarrow Z$ is
a non-isomorphism, then $ht = 0$.

Under (AR1) and (AR2),  (AR3) is equivalent to

(AR3') \ If $s: X\longrightarrow X'$ is not a section, then there exists
$s': Y\longrightarrow X'$ such that $s = s'f$.

Also, (AR3') is equivalent to

(AR4') \ If $X'$ is indecomposable and $s: X\longrightarrow X'$ is
a non-isomorphism, then $s\circ h[-1] = 0$.

\vskip10pt

In an Auslander-Reiten triangle $X {\longrightarrow} Y
{\longrightarrow} Z {\longrightarrow} X[1]$, the object $X$ is
uniquely determined by $Z$. Write $X = \tau_\mathcal A Z$. In
general $\tau_\mathcal A$ is {\em not} a functor. By definition
$\mathcal A$ has right Auslander-Reiten triangles if
there exists an
Auslander-Reiten triangle $X {\longrightarrow} Y {\longrightarrow}
Z {\longrightarrow} X[1]$ for any indecomposable $Z$; and $\mathcal A$ has
Auslander-Reiten triangles if $\mathcal{A}$ has right and
left Auslander-Reiten triangles.  We refer to [Hap1], [XZ] and [A] for the
Auslander-Reiten quiver of a triangulated category.

\vskip10pt

A Hom-finite $k$-category is {\em Krull-Schmidt} if the
endomorphism algebra of any indecomposable is local. In this case
any object is uniquely decomposed into a direct sum of
indecomposables, up to isomorphisms and up to the order of
indecomposable direct summands (Ringel [R], p.52).

Let $\mathcal{A}$ be a Hom-finite Krull-Schmidt triangulated
$k$-category. Theorem I.2.4 in [RV] says that $\mathcal{A}$ has a
right Serre functor $F$ if and if $\mathcal{A}$ has right
Auslander-Reiten triangles. In this case, $F$ coincides with
$[1]\circ \tau_\mathcal A$ on objects, up to isomorphisms.

\vskip10pt

\subsection{}Let $\mathcal{A}$ be a Hom-finite
triangulated $k$-category with Serre functor $F$. Denote by $[1]$
the shift functor of $\mathcal A$. Following Kontsevich [Ko],
$\mathcal{A}$ is {\em a Calabi-Yau category} if there is a natural
isomorphism $F \cong [d]$ of functors for some $d\in\Bbb Z$.

Denote by $o([1])$ the order of $[1]$. If $o([1]) = \infty$ then
the integer $d$ above is unique, and is called {\em the  CY dimension} of
$\mathcal{A}$; if $o([1])$ is finite then we call the minimal non-negative integer $d$
such that $F\cong [d]$ {\em the CY dimension} of
$\mathcal{A}$. Denote by $\operatorname{CYdim}(\mathcal A)$ the CY
dimension of $\mathcal A$.

\vskip10pt

For example, if $A$ is a symmetric algebra and $\mathcal P$ is the
category of projective modules, then the homotopy category
$K^b(\mathcal P)$ is of CY dimension $0$. Moreover, if $\mathcal
A$ is of CY dimension $d$, then $\operatorname {Ext}_\mathcal
A^i(X, Y)\cong D\circ \operatorname {Ext}_\mathcal A^{d-i}(Y, X),
\ X, Y\in\mathcal A, \ i\in\Bbb Z$, where $\operatorname
{Ext}_\mathcal A^i(X, Y): =\operatorname {Hom}_\mathcal A(X,
Y[i])$. Thus, if $A$ is a CY algebra ([B], [Gin2]), i.e. the
bounded derived category $D^b(A\mbox{-mod})$ is Calabi-Yau of CY
dimension $d$, then $\operatorname {gl.dim}(A\mbox{-mod})=d$ (see
[B]).

\vskip10pt

\subsection{} Let $\mathcal{A}$ and $\mathcal{B}$ be triangulated categories.
{\em A triangle functor} from $\mathcal{A}$ to $\mathcal{B}$ is a pair $(F, \eta^F)$, where
$F\colon \mathcal{A} \to \mathcal{B}$ is an additive functor and
$\eta^F: \ F\circ [1] \longrightarrow [1]\circ F$ is a natural
isomorphism, such that if $X \stackrel{f} {\longrightarrow} Y
\stackrel{g} {\longrightarrow} Z \stackrel{h} {\longrightarrow}
X[1]$ is a distinguished triangle of $\mathcal{A}$ then $FX
\stackrel{Ff} {\longrightarrow} FY \stackrel{Fg} {\longrightarrow}
FZ \stackrel{\eta^F_X \circ Fh} {\longrightarrow} (FX)[1]$ is a
distinguished triangle of $\mathcal{B}$. Triangle functors $(F, \
\eta^F)$ and $(G, \ \eta^G)$ are {\em natural isomorphic} if there
is a natural isomorphism $\xi: \ F\longrightarrow G$ such that the
following diagram commutes for any $A\in\mathcal{A}$
\[\xymatrix{
F(A[1])  \ar[rr]^{\eta^F_A} \ar[d]_{\xi_{A[1]}} && F(A)[1]
 \ar[d]^-{\xi_A[1]}\\
G(A[1]) \ar[rr]^{\eta^G_A} && G(A)[1].}\]

As Keller pointed out, the pair $([n], \ (-1)^n{\rm Id}_{[n+1]}):
\ \mathcal{A}\longrightarrow\mathcal{A}$ is a triangle functor for
$n\in\Bbb Z$. However, $([n], \ {\rm Id}_{[n+1]})$ may be {\em
not}.

We need the following important result. A nice proof given by Van
den Bergh is in the Appendix of [B].

\begin {lem} \ \ (Bondal-Kapranov {\em [BK]}; Van den
Bergh {\em [B]}) \ \ Let $F$ be a Serre functor of a Hom-finite
triangulated $k$-category $\mathcal{A}$. Then there exists a
natural isomorphism $\eta^F: \ F\circ [1] \longrightarrow [1]\circ
F$ such that $(F, \ \eta^F): \
\mathcal{A}\longrightarrow\mathcal{A}$ is a triangle functor.
\end{lem}

>From [Ke, 8.1] and [B, A.5.1] one has the following

\begin{prop} (Keller; Van den Bergh)\ Let $\mathcal{A}$ be a
Hom-finite triangulated $k$-category with Serre functor $F$. Then
$\mathcal{A}$ is a Calabi-Yau category if and only if there exists
a natural isomorphism $\eta^F: F\circ [1] \longrightarrow [1]\circ
F,$ such that $(F, \ \eta^F)$ is a triangle functor and $(F,
\eta^F) \longrightarrow ([d], (-1)^d\operatorname {Id}_{[d+1]})$\
is a natural isomorphism of triangle functors, for some integer
$d$. \end{prop}

\noindent {\bf Proof.} For convenience we justify the ``only if''
part. By assumption we have a natural isomorphism $\xi: F \cong
[d]$. Define $\eta^F_A \colon F(A[1]) \longrightarrow (FA)[1]$ for
$A\in\mathcal{A}$ by \ $\eta^F_A:=(-1)^d(\xi_A)^{-1}[1]\circ
\xi_{A[1]}.$\  Then $\xi_A[1] \circ \eta_A^F = (-1)^d \xi_{A[1]}$.
The naturality of $\eta^F: \ F\circ [1] \longrightarrow [1]\circ
F$ follows from the one of $\xi$. It remains to show that $(F,
\eta^F)\colon \mathcal{A}\to \mathcal{A}$ is a triangle functor.
Let $X \stackrel{f} {\longrightarrow} Y \stackrel{g}
{\longrightarrow}  Z \stackrel{h} {\longrightarrow} X[1]$ be a
distinguished triangle. Since \ $X[d]\stackrel{f[d]}
{\longrightarrow}Y[d]\stackrel{g[d]} {\longrightarrow} Z[d]
\stackrel{(-1)^d h[d]} {\longrightarrow}X[d+1]$\  \ is a
distinguished triangle, it suffices to prove that the following
diagram is commutative
\[\xymatrix{ F(X)\ar[r]^-{F(f)} \ar[d]_{\xi_X}^-{\wr} & F(Y)
\ar[r]^-{F(g)}\ar[d]_{\xi_Y}^-{\wr}& F(Z)\ar[rr]^-{\eta_X^F \circ
F(h)}\ar[d]_{\xi_Z}^-{\wr}&&
(F(X))[1]\ar[d]\ar[d]_{\xi_X[1]}^-{\wr}\\
X[d]\ar[r]^-{f[d]} & Y[d]\ar[r]^-{g[d]}& Z[d]\ar[rr]^-{(-1)^d
h[d]}&& X[d+1].}\] By the naturality of $\xi$ the first and the second square are
commutative. We also have
$$\xi_X[1]\circ \eta_X^F \circ F(h) = (-1)^d \xi_{X[1]} \circ F(h) = (-1)^d h[d]\circ \xi_Z.
\ \ \ \blacksquare $$

\vskip10pt

\subsection{} Let $A$ be a self-injective $k$-algebra,
$A$-mod the category of finite-dimensional left $A$-modules,
and $A\underline {\mbox{-mod}}$ the stable category of
$A$-mod modulo projective modules. Then the Nakayama functor
$\mathcal N: = D( A)\otimes_A-$, Heller's syzygy functor $\Omega$,
and the Auslander-Reiten translate $\tau
\cong \Omega^2\circ \mathcal N\cong\mathcal N\circ \Omega^2$
([ARS], p.126), and  are endo-equivalences of
$A\underline {\mbox{-mod}}$ ([ARS], Chap. IV). Note that
$A\underline {\mbox{-mod}}$ is a Hom-finite Krull-Schmidt
triangulated $k$-category with $[1] = \Omega^{-1}$ ([Hap1], p.16).
By the bi-naturality  of the Auslander-Reiten isomorphisms ([AR])
$$ \underline {\operatorname {Hom}}(X, Y) \cong D\circ\operatorname
{Ext}_A^1(Y, {\tau} X) \cong D\circ\underline {\operatorname
{Hom}}(Y, [1]\circ {\tau} X),$$  where $\underline {\operatorname
{Hom}}(X, Y): = \operatorname {Hom} _{A\underline
{\mbox{-mod}}}(X, Y)$, one gets the Serre functor
$F:=[1]\circ {\tau} \cong \Omega\circ \mathcal N$ of $A\underline {\mbox{-mod}}$.
It follows that $A\underline {\mbox{-mod}}$ is
Calabi-Yau if and only if $\mathcal N \cong
\Omega^{-(d+1)}$ for some $d$ ([Ke, 8.3]). In this case denote by
$\operatorname{CYdim}(A)$ the CY dimension of $A\underline
{\mbox{-mod}}$. Note that $\Omega, \ F, \ \mathcal N,
\ \tau$ are pairwise commutative
as functors of $A\underline {\mbox{-mod}}$. This follows from Lemma 2.1.
\vskip10pt

\subsection{} Let $A$ be a finite-dimensional
$k$-algebra. Recall that $A$ is a Nakayama algebra if any
indecomposable is uniserial, i.e. it has a unique composition
series ([ARS], p.197). In this case $A$ is representation-finite.
If $k$ is algebraically closed then any connected self-injective
Nakayama algebra is Morita equivalent to $\A(n, t),$ $n\ge 1, \
t\ge 2$ ([GR], p.243), which is defined below.

Let $\Bbb Z_n$ be the cyclic quiver with vertices indexed by the
cyclic group $\Bbb Z/n\Bbb Z$ of order $n$, and with arrows $a_i:
\ i \longrightarrow i+1, \ \forall \ i\in \Bbb Z/n\Bbb Z$. Let
$k\Bbb Z_n$ be the path algebra of the quiver $\Bbb Z_n$, $J$ the
ideal generated by all arrows, and $\A = \A(n, t): =k\Bbb Z_n/J^t$
with $t\ge 2$. Denote by $\gamma^l_i$ the path starting at vertex
$i$ and of length $l$, and $e_i: = \gamma^0_i$. We write the
conjunction of paths from right to left. Then $\{\gamma^l_i \ | \
0\le i\le n-1, \ 0\le l\le t-1\}$ is a basis of $\A$; while
$\{P(i): = \A e_i\ | \ 0\le i\le n-1\}$ is the set of pairwise
non-isomorphic indecomposable projective modules, and $\{I(i): =
D(e_{i}\A)\ | \ 0\le i\le n-1\}$ is the set of pairwise
non-isomorphic indecomposable injective modules, with $P(i) \cong
I(i+t-1)$. Note that $\A$ is a Frobenius algebra, and $\A$ is
symmetric if and only if $n\mid (t-1)$. Write $S(i): =
P(i)/\operatorname{rad}P(i)$, and $S^l_i: =
\A\gamma_{i+l-t}^{t-l}$. Then $S_i^l$ is the indecomposable with
top $S(i)$ and the Loewy length $l$, and $\{ S^l_i \ | \ 0\le i\le
n-1, \ 1\le l\le t \}$ is the set of pairwise non-isomorphic
indecomposable modules, with $S^{t}_i = P(i)$ and
$\operatorname{soc}(S^l_i) = S(i+l-1)$. For the Auslander-Reiten
quiver of $\A$ see [GR], Section 2, and [ARS], p.197. In
particular, the stable Auslander-Reiten quiver of $\A$ is $\Bbb Z
A_{t-1}/\langle\tau^n\rangle.$

\vskip10pt

\section{\bf Indecomposable Calabi-Yau objects}

The purpose of this section is to introduce the Calabi-Yau objects
and to give the relation between indecomposable
Calabi-Yau objects and Auslander-Reiten triangles.

\vskip10pt

\subsection{} Let $\mathcal{A}$ be a Hom-finite
triangulated $k$-category. A non-zero object $X$ is called \emph
{a Calabi-Yau object} if there exists  a natural isomorphism
\begin{align}\operatorname{Hom}_{\mathcal A}(X, -)\cong D\circ \operatorname{Hom}_{\mathcal A}(-,
X[d])\end {align} for some integer $d$.

By Yoneda Lemma, such a $d$ is unique up to a multiple of the
relative order $o([1]_X)$ of $[1]$ respect to $X$. Recall that
$o([1]_X)$ is the minimal positive integer such that
$X[o([1]_X)]\cong X$, otherwise $o([1]_X) = \infty$. If
$o([1]_X)=\infty$ then $d$ in $(3.1)$ is unique and is called {\em
the CY dimension} of $X$. If $o([1]_X)$ is finite then the minimal
non-negative integer $d$ in $(3.1)$ is called {\em the CY
dimension} of $X$. We denote $\operatorname{CYdim}(X)$ the CY
dimension. Thus, if $o([1]) < \infty$ then $o([1]_X)\mid o([1])$
and $0\le \operatorname{CYdim}(X) < o([1]_X).$

Let $A$ be a finite-dimensional self-injective algebra. An
$A$-module $M$ without projective direct summands is called {\em a
Calabi-Yau module} of CY dimension $d$, if it is a Calabi-Yau
object of $A\underline{\mbox{-mod}}$ with $\operatorname{CYdim}(M) = d$.

\vskip10pt

Note that $\operatorname{CYdim}(X)$ is usually not easy to
determine. In case $(3.1)$ holds for some $d$, we say that $X$ is a {\em  $d$-th CY object}.
Of course, if $o([1]_X) <
\infty$ then $o([1]_X)\mid (d - \operatorname{CYdim}(X))\ge 0.$

\vskip10pt

If $\mathcal A$ has right Serre functor $F$, then by Yoneda Lemma
a non-zero object $X$ is a $d$-th CY object if and only if $F(X)\cong X[d]$, or equivalently,
$F(X)[-d]\cong X$. Thus, a non-zero $A$-module $M$
without projective direct summands is a $d$-th CY module
if and only if $\mathcal N(M)\cong
\Omega^{-(d+1)}(M)$ in $A\underline{\mbox{-mod}}$ (in fact, this
isomorphism can be taken in $A$-mod).

\vskip10pt

\subsection{} We have the following basic property.

\vskip10pt

\begin {prop} \ $(i)$ \ The Calabi-Yau property for a category or an object,
is invariant under triangle-equivalences.

\vskip10pt

$(ii)$ \ The Calabi-Yau property for a module is {\em ``usually"}
invariant under stable equivalences between self-injective
algebras.  Precisely, let $A$ and $B$ be self-injective
algebras, $G: A\underline {\mbox{-mod}}\longrightarrow B\underline
{\mbox{-mod}}$ a stable equivalence, and $X$ a CY $A$-module of
dimension $d$. If $A\ncong \A(n, 2),$ or if $A$ and $B$ are
symmetric algebras, then $G(X)$ is a CY $B$-module of dimension
$d$.
\end{prop}

\noindent{\bf Proof.} \ $(i)$ \ Let $\mathcal A$ be a Calabi-Yau
category with $F_\mathcal A\cong [d]$, where $F_\mathcal A$ is the
Serre functor. Clearly $F_\mathcal B: = G \circ
F_\mathcal A\circ G^{-1}$ is a Serre functor of $\mathcal B$ (if
$\mathcal B$ has already one, then it is natural
isomorphic to $F_\mathcal B$). By the natural
isomorphism $(\xi^G)^d: G\circ [d] \longrightarrow [d]\circ G$,
which is the composition $G\circ [d] \longrightarrow [1]\circ
G\circ [d-1] \longrightarrow \cdots \longrightarrow [d]\circ G$
(A.2 in [B]), we see that $\mathcal B$ is a Calabi-Yau
category with $F_\mathcal B\cong [d].$ If $X$ is a calabi-Yau object
with a natural isomorphism $\eta$ as in
$(3.1)$, then we have natural isomorphism \
$\operatorname{Hom}_\mathcal B(-, (\xi^G)^d_X)\circ G\circ
\eta\circ G^{-1}: \ \operatorname{Hom}_{\mathcal B}(G(X), -)\cong
D\circ \operatorname{Hom}_{\mathcal B}(-, (GX)[d]),$ \ which
implies that $G(X)$ is a Calabi-Yau object of $\mathcal B$.

$(ii)$ \ Recall that an equivalence $G: A\underline
{\mbox{-mod}}\longrightarrow B\underline {\mbox{-mod}}$ of
categories is called a stable equivalence. Note that in general
$G$ is not induced by an exact functor (cf. [ARS], p.339), hence
$G$ may be not a triangle-equivalence (cf. [Hap1], Lemma 2.7,
p.22. Note that the converse of Lemma 2.7 is also true). One may
assume that $A$ is connected. If $A\ncong \A(n, 2),$ or if $A$ and
$B$ are symmetric algebras, then by Corollary 1.7 and Prop. 1.12
in [ARS], p.344, we know that $G$ commutes with $\tau$ and
$\Omega$ on modules, hence we have isomorphism
\begin{align*}\Omega_B^{-1}\circ \tau_B (G(X)) \cong
G(\Omega_A^{-1}\circ \tau_A) (X))  \cong G(\Omega_A^{-d}(X))\cong
\Omega_B^{-d}(G(X)),\end{align*} which implies that $G(X)$ is a
Calabi-Yau $B$-module of CY dimension $d$. \hfill $\blacksquare$

\vskip10pt

It seems that the Calabi-Yau property for the stable category
is also invariant under stable equivalence $G$
between self-injective algebras. However, this need natural isomorphiams between
$G$ and $\tau$, and $G$ and $\Omega$, which are not clear to us.

\vskip10pt

\subsection{} The main result of this section
is as follows.

\vskip10pt

\begin{thm} Let $\mathcal{A}$ be a Hom-finite Krull-Schmidt triangulated
$k$-category, and $X$ an indecomposable object of $\mathcal A$.
Then $X$ is a $d$-th CY object if and
only if there exists an Auslander-Reiten triangle of the form
\begin{align}X[d-1]\stackrel f \longrightarrow Y
\stackrel g \longrightarrow X\stackrel h\longrightarrow
X[d].\end{align}

Moreover, $Y$ is also a $d$-th CY object.
\end{thm}

\vskip10pt
\subsection{}The proof of the first part of Theorem 3.2
follows an argument of Reiten and Van den Bergh in [RV]. For the convenience
we include a complete proof.

\vskip10pt

\begin{lem} \ Let $\mathcal{A}$ be a Hom-finite Krull-Schmidt triangulated
$k$-category, and $X$ a non-zero object of $\mathcal A$. Then $X$
is a $d$-th CY object if and only if for
any indecomposable $Z$ there exists a non-degenerate bilinear form
\begin{align}(-,-)_{Z}: \
\operatorname{Hom}_{\mathcal A}(X, Z)\times
\operatorname{Hom}_{\mathcal A}(Z, X[d]) \longrightarrow
k\end{align} such that for any \ $u\in
\operatorname{Hom}_{\mathcal A}(X, Z), \ \ v\in
\operatorname{Hom}_{\mathcal A}(Z, W)$, \ and \ $w\in
\operatorname{Hom}_{\mathcal A}(W, X[d]),$ there holds
\begin{align}(u, wv)_Z=(vu, w)_W.\end{align}

\end{lem}
\noindent{\bf Proof.} \ If $X$ is a $d$-th CY object, then  we
have $\eta_Z: \ \operatorname{Hom}_{\mathcal A}(X, Z)\cong D\circ
\operatorname{Hom}_{\mathcal A}(Z, X[d]), \ \forall \ Z\in\mathcal
A$,  which are natural in $Z$. Each isomorphism $\eta_Z$ induces a
non-degenerate bilinear form $(-,-)_{Z}$ in $(3.3)$ by $(u, z)_Z:
= \eta_Z(u)(z),$ and $(3.4)$ follows from the naturality of
$\eta_Z$ in $Z$. Conversely, if we have $(3.3)$ and $(3.4)$ for
any indecomposable $Z$, then we have isomorphism $\eta_Z: \
\operatorname{Hom}_{\mathcal A}(X, Z)\cong D\circ
\operatorname{Hom}_{\mathcal A}(Z, X[d])$ given by $\eta_Z(u)(z):
= (u, z)_Z, \ \forall \ z\in \operatorname{Hom}_{\mathcal A}(Z,
X[d])$.  By $(3.4)$  $\eta_Z$ are natural in $Z$. Since $\mathcal
A$ is Krull-Schmidt, it follows that we have isomorphisms $\eta_Z$
for any $Z$ which are natural in $Z$. This means that $X$ is a
$d$-th CY object. \hfill $\blacksquare$

\vskip10pt

The following Lemma in [RV] will be used.

\vskip10pt

\begin{lem} \ \ ({\em [RV]}, Sublemma I.2.3)\  Let $\mathcal{A}$ be
a Hom-finite Krull-Schmidt triangulated $k$-category,
$\tau_{\mathcal A}(X)  \longrightarrow Y
  \longrightarrow X\stackrel h\longrightarrow
\tau_\mathcal A(X)[1]$ an Auslander-Reiten triangle of $\mathcal
A$, and $Z$ an indecomposable in $\mathcal A$. Then

$(i)$ \ For any non-zero $z\in \operatorname{Hom}_{\mathcal A}(Z,
\tau_{\mathcal A}(X)[1])$ there exists $u\in
\operatorname{Hom}_{\mathcal A}(X, Z)$ such that $zu = h$.

$(ii)$ For any non-zero $u\in \operatorname{Hom}_{\mathcal A}(X,
Z)$ there exists $z\in \operatorname{Hom}_{\mathcal A}(Z,
 \tau_{\mathcal A}(X)[1])$ such that $zu = h$.
\end{lem}

\vskip10pt

In a Hom-finite Krull-Schmidt triangulated $k$-category without
Serre functor (e.g., by [Hap2] and [RV] if
$\operatorname{gl.dim}(A) = \infty$ then $D^b(A\mbox{-mod})$ has
no Serre functor), one may use {\em the generalized Serre functor}
introduced by Chen [Ch].

\vskip10pt

\begin{lem} (Chen [Ch]) Let $\mathcal A$ be a Hom-finite Krull-Schmidt triangulated
$k$-category. Consider the full subcategories of $\mathcal A$
given by $$\mathcal{\mathcal A}_r:=\{X \in \mathcal{A}\; |\;
D\circ\operatorname{Hom}_{\mathcal A}(X, -) \mbox{ is
representable} \}$$ and $$\mathcal{\mathcal A}_l:=\{X \in
\mathcal{A}\;|\; D\circ \operatorname{Hom}_\mathcal A(-, X) \mbox{
is representable}\}.$$ Then both $\mathcal{A}_r$ and
$\mathcal{A}_l$ are thick triangulated subcategories of $\mathcal
A$. Moreover, one has

$(i)$ \ There is a unique $k$-functor $S: \mathcal{A}_r
\longrightarrow \mathcal{A}_l$ which is an equivalence,  such that
there are natural isomorphisms
\begin{align}\operatorname{Hom}_\mathcal A(X, -) \simeq
D\circ\operatorname{Hom}_\mathcal A(-, S(X)), \ \ \forall \ X\in
\mathcal {A}_r,\end{align} which are natural in $X$. $S$ is called
the generalized Serre functor, with range $\mathcal A_r$ and
domain $\mathcal A_l$.

$(ii)$ \ There exists a natural isomorphism $\eta^S: S\circ
[1]\longrightarrow [1]\circ S$ such that the pair $(S, \eta^S):
\mathcal{A}_r \longrightarrow \mathcal{A}_l$ is an
triangle-equivalence.
\end{lem}

In this terminology, a non-zero object $X$ is a $d$-th CY object
if and only if $X\in\mathcal A_r$ and $S(X)\cong X[d]$, by $(3.5)$
and Yoneda Lemma.

\vskip10pt

\subsection{} {\bf Proof of Theorem 3.2.}\quad Let
$X$ be an indecomposable $d$-th CY object. By Lemma 3.3 we have a
non-degenerate bilinear from $(-,-)_X: \
\operatorname{Hom}_{\mathcal A}(X, X)\times
\operatorname{Hom}_{\mathcal A}(X, X[d]) \longrightarrow k.$ It
follows that there exists $0\ne h\in \operatorname{Hom}_{\mathcal
A}(X, X[d])$ such that $(\operatorname{radHom}_{\mathcal A}(X, X),
h)_X=0.$ Embedding $h$ into a distinguished triangle as in
$(3.2)$. We claim that it is an Auslander-Reiten triangle. For
this it remains to prove (AR4) in 2.1. Let $X'$ be indecomposable
and $t: X'\longrightarrow X$ a non-isomorphism. Then by $(3.4)$
for any $u\in \operatorname{Hom}_{\mathcal A}(X, X')$ we have $(u,
ht)_{X'} = (tu, h)_X = 0.$ Since $(-,-)_{X'}$ is non-degenerate,
it follows that $ht=0$.

Conversely, let $(3.2)$ be an Auslander-Reiten triangle. In order
to prove that $X$ is a $d$-th CY object, by Lemma 3.3 it suffices
to prove that for any indecomposable $Z$ there exists a
non-degenerate bilinear form $(-,-)_{Z}$ as in $(3.3)$ satisfying
$(3.4)$. For this, choose an arbitrary linear function
$\operatorname {tr}\in D\circ \operatorname{Hom}_{\mathcal A}(X,
X[d])$ such that $\operatorname {tr}(h)\ne 0$, and define $(u,
z)_Z = \operatorname{tr}(zu).$ Then $(3.4)$ is automatically
satisfied. It remains to prove that $(-, -)_Z$ is non-degenerate.
In fact, for any $0\ne z\in \operatorname{Hom}_{\mathcal A}(Z,
X[d])$, by Lemma 3.4 there exists $u\in
\operatorname{Hom}_{\mathcal A}(X, Z)$ such that $zu = h$. So $(u,
z)_Z = \operatorname{tr}(zu) = \operatorname{tr}(h)\ne 0$.
Similarly, for any $0\ne u\in \operatorname{Hom}_{\mathcal A}(X,
Z)$ we have $z\in \operatorname{Hom}_{\mathcal A}(Z, X[d])$ such
that $(u, z)_Z \ne 0$. This proves the non-degenerateness of $(-,
-)_Z$.

\vskip10pt

Now we prove that $Y$ in $(3.2)$ is also a $d$-th CY object. We make use of the
generalized Serre functor in [Ch]
(For the reader prefer Serre functor, one can assume the existence, and use Lemma 2.1).
Since $X, \ X[d-1]\in\mathcal A_r$, it follows from Lemma 3.5 that
$Y\in\mathcal A_r$. Applying the
generalized Serre functor $(S, \eta^S)$ to $(3.2)$ we get the
distinguished triangle (by Lemma 3.5$(ii)$)
$$S(X[d-1])\stackrel {S(f)}\longrightarrow S(Y)\stackrel {S(g)}\longrightarrow
S(X)\stackrel {\eta^S_{X[d-1]}\circ S(h)}\longrightarrow
S(X[d-1])[1].\eqno(*)$$ Also, we have the Auslander-Reiten
triangle
$$X[2d-1]\stackrel {f[d]}
{\longrightarrow}Y[d] \stackrel {g[d]}
{\longrightarrow}X[d]\stackrel {(-1)^d h[d]}{\longrightarrow}
X[2d].$$ Since $X$ is a $d$-th CY object it
follows that we have an isomorphism $w: S(X)\longrightarrow X[d]$.
Note that $S(h)\ne 0$ means that $S(g)$ is not a retraction
([Hap1], p.7). Thus, by (AR3) there exists $v: S(Y)\longrightarrow
Y[d]$ such that $w\circ S(g) = g[d]\circ v$. By the definition of a triangulated category
we get $u:
S(X[d-1])\longrightarrow X[2d-1]$ such that the following diagram
is commutative
\[\xymatrix{ S(X[d-1])\ar[r]^-{S(f)}
\ar[d]_{u} & S(Y) \ar[r]^-{S(g)}\ar[d]_{v}&
S(X)\ar[rr]^-{\eta_{X[d-1]}^S \circ S(h)}\ar[d]_{w}&&
S(X[d-1])[1]\ar[d]\ar[d]_{u[1]}\\
X[2d-1]\ar[r]^-{f[d]} & Y[d]\ar[r]^-{g[d]}& X[d]\ar[rr]^-{(-1)^d
h[d]}&& X[2d].}\eqno(**)\] We claim that $u: S(X[d-1])\longrightarrow
X[d]$ is an isomorphism, hence by the property of a triangulated
category we know that $v: S(Y)\longrightarrow Y[d]$ is also an
isomorphism, i.e. $Y$ is a $d$-th CY object.

Otherwise $u: S(X[d-1])\longrightarrow X[2d-1]$ is not an
isomorphism. Note that $S$ is {\em only} defined on $\mathcal A_r$, it follows that
we do not know if $(*)$ is an Auslander-Reiten triangle.

Since $X[d-1]$ is a $d$-th CY object, it
follows that we have an isomorphism $\alpha:
X[2d-1]\longrightarrow S(X[d-1])$, hence $\alpha\circ u\in
\operatorname{Hom}_\mathcal A(S(X[d-1]), S(X[d-1]))$. So we have
$u'\in \operatorname{Hom}_\mathcal A(X[d-1], X[d-1])$ such that $
\alpha\circ u = S(u')$, and $u'$ is also a non-isomorphism. Since
$X[d-1]\stackrel f \longrightarrow Y \stackrel g \longrightarrow
X\stackrel h\longrightarrow X[d]$ is an Auslander-Reiten triangle
it follows from (AR4') that $u'\circ h[-1] = 0$, or equivalently,
$S(u')\circ S(h[-1]) = 0$ (note that $h[-1]\in \mathcal A_r$).
Thus we have
$$u[1] \circ S(h[-1])[1]\circ \eta^S_{X[-1]} = 0$$
where $\eta^S_{X[-1]}: S(X)\longrightarrow S(X[-1])[1]$ is an
isomorphism. By the naturality of $\eta^S$ we have the
commutative diagram
\[\xymatrix{
S(X) \ar[rr]^{\eta_{X[-1]}^S} \ar[d]_{S(h)}&&
S(X[-1])[1] \ar[d]^-{S(h[-1])[1]}\\
S(X[d]) \ar[rr]^{\eta_{X[d-1]}^S} && S(X[d-1])[1].}\] It follows
that we have $u[1]\circ \eta^S_{X[d-1]}\circ S(h) = 0,$ and hence
by the commutative diagram $(**)$ we get a contradiction \ $(-1)^dh[d]\circ w = u[1]\circ
\eta^S_{X[d-1]}\circ S(h) = 0.$ This completes the proof. \hfill $\blacksquare$

\vskip10pt

\subsection{} {\bf Remark 3.6.} Let $\mathcal A$ be a Hom-finite Krull-Schmidt triangulated
$k$-category. If every indecomposable $X$ in $\mathcal A$ is a
$d_X$-th CY object, then $\mathcal A$ has a Serre functor $F$ with
$F(X)\cong X[d_X]$.

In fact, by Theorem 3.2 $\mathcal A$ has right and left Auslander-Reiten triangles,
and then by Theorem I.2.4 in [RV] $\mathcal A$ has Serre functor
$F$. By Prop.I.2.3 in [RV] and $(3.2)$ we have $F(X)\cong
\tau_\mathcal A(X)[1] = X[d_X-1][1] = X[d_X]$ for any indecomposable
$X$.

\vskip10pt

However, even if all indecomposables are $d$-th CY objects with
the same $d$, we do not know whether $\mathcal A$ is a Calabi-Yau
category, although $F$ and $[d]$ coincide on objects. The examples
we know have a positive answer to this question. \hfill
$\blacksquare$

\vskip10pt

\section{\bf Minimal Calabi-Yau objects}

The purpose of this section is to describe all the Calabi-Yau objects of a
Hom-finite Krull-Schmidt triangulated $k$-category with a Serre functor.

\subsection{} Let $\mathcal{A}$ be a
Hom-finite Krull-Schmidt triangulated $k$-category. A $d$-th CY
object $X$ is said to be {\em minimal} if
any {\em proper} direct summand of $X$ is {\em not} a $d$-th CY
object.

\vskip10pt

\begin{lem} \ Let $\mathcal{A}$ be a Hom-finite Krull-Schmidt triangulated
$k$-category with right Serre functor $F$. Then a non-zero object
$X$ is a minimal $d$-th CY object if and
only if the following  are satisfied:

\vskip10pt

1. The indecomposable direct summands of $X$ can be ordered as $X
= X_1\oplus \cdots \oplus X_r$ such that
\begin{align} F(X_1) \cong X_2[d],
\ F(X_2) \cong X_3[d],  \ \cdots, \ F(X_{r-1}) \cong X_r[d],  \
F(X_r)\cong X_1[d].\end{align}

We call the cyclic order arising from this property {\em a canonical order} of $X$
(with respect to $F$ and $[d]$).

2. \ $X$ is multiplicity-free, i.e. its indecomposable direct
summands are pairwise non-isomorphic.
\end{lem}

\noindent {\bf Proof.}\quad In the following we often use that a
non-zero object $X$ is a $d$-th CY object if
and only if $F(X)\cong X[d]$.

Let $X= X_1\oplus \cdots \oplus X_r$ be a minimal $d$-th CY
object, with each $X_i$ indecomposable and
$r\ge 2$. Then $F(X_1) \oplus\cdots\oplus F(X_r)\cong
X_1[d]\oplus\cdots\oplus X_r[d]$. Since $\mathcal{A}$ is
Krull-Schmidt, it follows that there exists a permutation $\sigma$
of $1, \cdots, r$, such that $F(X_i)\cong X_{\sigma(i)}[d]$ for
each $i$. Write $\sigma$ as a product of disjoint cyclic
permutations. Since $X$ is minimal, it follows that $\sigma$ has
to be a cyclic permutation of length $r$. By reordering the
indecomposable direct summands of $X$, one may assume that $\sigma
= (12\cdots r)$. Thus, $X$ satisfies the condition 1.

Now, we consider a canonical order $X= X_1\oplus \cdots \oplus
X_r$. If $X_i\cong X_j$ for some $1\le i < j\le r$, then
$F(X_{j-1})\cong X_j[d]\cong X_i[d]$, it follows that $X_i\oplus
\cdots\oplus X_{j-1}$ is already a $d$-th CY object,
which contradicts the minimality of $X$. This
proves that $X$ is multiplicity-free.

\vskip10pt Conversely, assume that a multiplicity-free object $X=
X_1\oplus \cdots \oplus X_r$ is in a canonical order. By $(4.1)$
we have $F(X)\cong X[d]$. So $X$ is a $d$-th CY object.
It remains to show the minimality.
If not, then there exists a proper direct summand $X_{i_1}\oplus
\cdots\oplus X_{i_t}$ of $X$ which is a minimal $d$-th CY
object, so $1\le t< r$. By what we have proved above
we may assume that this is a canonical order. Then
$$F(X_{i_1})\cong X_{i_2}[d],\  \cdots, \ F(X_{i_{t-1}})\cong
X_{i_t}[d], \ F(X_{i_t})\cong X_{i_1}[d].$$ While $X= X_1\oplus
\cdots \oplus X_r$ is also in a canonical order, it follows that
(note that we work on indices modulo $r$, e.g. if $i_1 =
r$ then $i_1+1$ is understood to be $1$)
$$F(X_{i_1})\cong X_{i_1+1}[d], \ \cdots, \ F(X_{i_{t-1}})\cong
X_{i_{t-1}+1}[d], \ F(X_{i_t})\cong X_{i_t+1}[d].$$ Since $X$ is
multiplicity-free, it follows that (considering indices modulo
$r$)
$$i_2 = i_1+1, \ \cdots, \ i_t = i_{t-1}+1, \ i_1 = i_t + 1,$$
hence $i_1 = i_1+t$, which means $r\mid t$. This is impossible
since $1\le t< r.$ \hfill $\blacksquare$

\vskip10pt

\subsection{} Let $\mathcal{A}$ be a Hom-finite Krull-Schmidt triangulated
$k$-category with Serre functor $F$. For each $d\in\Bbb Z$,
consider the triangle-equivalence $G:=[-d]\circ F\cong
F\circ [-d]:  \mathcal A\longrightarrow \mathcal A$. For each
indecomposable $M\in\mathcal A$, denote by $o(G_M)$ the relative
order of $G$ respect to $M$, that is, $r:=o(G_M)$ is the
minimal positive integer such that $G^r(M)\cong M$, otherwise $o(G_M) =
\infty$. Denote by
$\operatorname{Aut}(\mathcal A)$ the group of the
triangle-equivalences of $\mathcal A$, and by $\langle G\rangle$
the cyclic subgroup of $\operatorname{Aut}(\mathcal A)$ generated
by $G$. Then $\langle G\rangle$ acts naturally on
$\operatorname{Ind}(\mathcal A)$, the set of the isoclasses of
indecomposables of $\mathcal A$. Denote by $\mathcal O_M$ the
$G$-orbit of an indecomposable $M$. Then $|\mathcal O_M| = o(G_M)$.
If $|\mathcal O_M| < \infty$ then the set $\mathcal O_M$ is a finite
$G$-orbit.

\vskip10pt

Denote by $\operatorname{Fin}\mathcal{O}(\mathcal A, d)$ the set
of all the finite $G$-orbits of $\operatorname{Ind}(\mathcal A)$,
and by $\operatorname{MinCY}(\mathcal A, d)$ the set of isoclasses
of minimal $d$-th CY objects. We
have the following

\vskip10pt

\begin {thm} \ \ Let $\mathcal{A}$ be a Hom-finite Krull-Schmidt triangulated
$k$-category with Serre functor $F$. Then

\vskip10pt

$(i)$ \ Every $d$-th CY object is a direct sum of finitely many
minimal $d$-th CY objects.

\vskip10pt

$(ii)$ \ With the notations above, for each $d\in\Bbb Z$ the map
\begin{align}\mathcal O_M \ \mapsto \
\bigoplus\limits_{X\in \mathcal O_M} X = M\oplus G(M)\oplus\cdots
\oplus G^{o(G_M)-1}(M)\end{align} \noindent gives a one-to-one
correspondence between the sets
$\operatorname{Fin}\mathcal{O}(\mathcal A, d)$ and
$\operatorname{MinCY}(\mathcal A, d)$, where $G: = [-d]\circ F$.

Thus, a minimal $d$-th CY object is
exactly the direct sum of all the objects in a finite $G$-orbit
of $\operatorname{Ind}(\mathcal A)$.

\vskip10pt

$(iii)$ Non-isomorphic minimal $d$-th CY objects are disjoint, i.e. they have no isomorphic
indecomposable direct summands. \end {thm}

\noindent {\bf Proof.}\quad Let $X$ be a $d$-th CY object.
If $X$ is not minimal, then $X = Y\oplus Z$
with $Y\ne 0\ne Z$ such that $F(Y)\oplus F(Z) \cong Y[d]\oplus
Z[d]$ and $F(Y)\cong Y[d]$. Since $\mathcal{A}$ is Krull-Schmidt,
it follows that $F(Z) \cong Z[d]$, i.e. $Z$ is also a $d$-th CY
object. Then $(i)$ follows by induction.

Thanks to Lemma 2.1, $(4.1)$ becomes $X_i = G^{i-1}(X_1), \ 1\le
i\le r.$ Then $(ii)$ is a reformulation of Lemma 4.1. Moreover
$(iii)$ follows from $(ii)$. \hfill $\blacksquare$

\vskip10pt

\begin {cor} \ \ Let $A$ be a finite-dimensional self-injective
algebra. Then $X$ is a minimal $d$-th CY module
if and only if $X$ is of the form $X \cong
\bigoplus \limits_{0\le i\le r-1}G^i(M)$, where $M$ is an
indecomposable non-projective $A$-module with $r: = o(G_M)<\infty$,
and $G: = \Omega^{d+1}\circ \mathcal N$. \end {cor}

\noindent {\bf Proof.}\quad Note that in this case $[-d]\circ F =
\Omega^{d+1}\circ \mathcal N$. \hfill $\blacksquare$

\vskip10pt

\subsection{} As an example, we describe all the
Calabi-Yau objects in $D^b(kQ\mbox{-mod})$, the bounded derived
category of $kQ\mbox{-mod}$, where $Q$ is a finite quiver without
oriented cycles.

Note that indecomposable objects of $D^b(kQ\mbox{-mod})$ are
exactly stalk complexes of indecomposable $kQ$-modules. The category
$D^b(kQ\mbox{-mod})$ has Serre functor $F =[1]\circ \tau_D$, where
$\tau_D$ is the Auslander-Reiten translation of
$D^b(kQ\mbox{-mod})$. Recall that $\tau_D$ is given by
$$\tau_D(M) =
\begin{cases}\tau(M), & \ \mbox{if} \ M \ \mbox{is an
indecomposable non-projective}; \\ I[-1], & \ \mbox{if}  \ P \
\mbox{is an indecomposable projective},
\end{cases}$$ where $I$ is the indecomposable injective with \ $\operatorname{soc}(I) =
P/\operatorname{rad} P$, and $\tau$ is the Auslander-Reiten
translation of $kQ\mbox{-mod}$ ([Hap1], p.51). Note that
$D^b(kQ\mbox{-mod})$ is {\em not} a Calabi-Yau category except
that $Q$ is the trivial quiver with one vertex and no arrows.
However, the cluster category $\mathcal C_{kQ}$
introduced in [BMRRT], which is the orbit category of
$D^b(kQ\mbox{-mod})$ respect to the functor $\tau_D^{-1}\circ
[1]$, is a Calabi-Yau category of CY dimension $2$ ([BMRRT]).

Let $M$ be a minimal Calabi-Yau object of $D^b(kQ\mbox{-mod})$ of
CY dimension $d$ (in this case $o([1]_M) = \infty$, hence $d$
is unique). By shifts we may assume that $M$ is a $kQ$-module. By
$F(M)=\tau_D(M[1]) = \tau_D(M)[1] = M[d]$ we see $d = 1$ or $0$.
Note that $kQ$ admits an indecomposable projective-injective
module if and only if $Q$ is of type $A_n$ with the linear
orientation. However, in this case the unique indecomposable
projective-injective module $P = I$ does not satisfy the relation
$\operatorname{soc}(I) = P/\operatorname{rad} P$. It follows that
$d\ne 0$. Thus $d=1$ and $\tau(M) = M$. Consequently, $Q$ is an
affine quiver and $M$ is a $\tau$-periodic (regular) module of
period $1$. All such modules are well-known, by the
classification of representations of affine quivers (see
Dlab-Ringel [DR]). Thus we have

\vskip10pt

\begin{prop} \ \ $D^b(kQ\mbox{-mod})$ admits a Calabi-Yau object $M$ if and only if
$Q$ is an affine quiver. In this case, $M$ is minimal if and
only if $M$ is an indecomposable in a homogeneous tube of
the Auslander-Reiten quiver of $kQ$-mod, or $M$ is the direct sum of all the indecomposables of
same quasi-length in a non-homogeneous tube of the Auslander-Reiten quiver of $kQ$-mod, up to
shifts. Moreover, all such $M$'s have CY dimension $1$.
\end{prop}

\vskip10pt

\section{\bf Calabi-Yau modules of self-injective Nakayama algebras}

This purpose of this section is to classify all the $d$-th CY
modules of self-injective Nakayama algebras $\A(n, t), \ n\ge 1, \
t\ge 2$, where $d$ is any given integer. By Theorem 4.2$(i)$ it
suffices to consider the minimal $d$-th CY $\A$-modules. By
Corollary 4.3 this reduces to computing the relative order $o(G_M)$
of $G: = \Omega^{d+1}\circ \mathcal N$ respect to any
indecomposable $\A(n, t)$-module $M=S_i^l, \ 1\le l\le t-1$, for
any integer $d$.

\vskip10pt

\subsection{}  Recall that $\A(n, t)$ is the quotient of the path algebra of the
cyclic quiver with $n$ vertices by the truncated ideal $J^t$,
where $J$ is the two-sided ideal generated by the arrows. From now
on we write $\A$ instead of $\A(n, t)$.

We keep the notations introduced in 2.5. Note that the
indecomposable module $S^l_i$ \ ($i\in\Bbb Z/n\Bbb Z,$ \ $1\le
l\le t-1$) \ has a natural $k$-basis, consisting of all the paths
of quiver $\Bbb Z_n$ starting at the vertex $i+l-t$ and of lengths
at least $t-l$:
$$\gamma^{t-l}_{i+l-t}, \ \gamma^{t-l+1}_{i+l-t}, \ \cdots, \ \gamma^{t-1}_{i+l-t}.$$
For $1\le l\le t-1$ denote by $\sigma_{i}^{l}:
S_{i}^{l}\longrightarrow S_{i-1}^{l+1}$ the inclusion by embedding
the basis above;  and for $2\le l\le t$ denote by $p_{i}^{l}:
S_{i}^{l}\longrightarrow S_{i}^{l-1}$ the $\A$-epimorphism given
by the right multiplication by arrow $a_{i+l-t-1}$. These
$\sigma_i^l$'s and $p_i^l$'s are all the irreducible maps of
$\A$-mod, up to scalars.

\vskip10pt

We need the explicit actions of functors $\mathcal N$ and
$\Omega^{-(d+1)}$ of $\A\underline {\mbox{-mod}}$. By the exact
sequences  $0 \longrightarrow S_{i}^l \longrightarrow
I(i+l-1)\longrightarrow S_{i+l-t}^{t-l}\longrightarrow 0$ (with
the canonical maps), via the  basis above one has the actions of
functor $\Omega^{-1}$ for $i\in\Bbb Z/n\Bbb Z, \ \ 1\le l \le
t-1$:
\begin{align*}\Omega^{-1}(S^l_i) =
S_{i+l-t}^{t-l}, \ \ \ \ \Omega^{-1}(\sigma^l_i) = p_{i+l-t}^{t-l},
\ \ \ \ \Omega^{-1}(p^l_i) = \sigma_{i+l-t}^{t-l}.\end{align*} By
induction one has in $\A\underline {\mbox{-mod}}$ for any integer
$m$ (even negative):
\begin{align} \Omega^{-(2m-1)}(S_i^l) = S_{i+l-mt}^{t-l}; \ \ \ \
\Omega^{-2m}(S_i^l) = S_{i-mt}^{l};\end{align}
and $$\Omega^{-(2m-1)}(\sigma_i^l) = p_{i+l-mt}^{t-l}, \ \ \ \
\Omega^{-(2m-1)}(p_i^l) = \sigma_{i+l-mt}^{t-l};\eqno(*)$$ and
$$\Omega^{-2m}(\sigma_i^l) = \sigma_{i-mt}^{l}, \ \ \ \
\Omega^{-2m}(p_i^l) = p_{i-mt}^{l}.\eqno(**)$$ In particular, we
have
$$o([1]) = \begin{cases} n, & \ t = 2;\\ 2m, & \ t\ge 3,\end{cases}\eqno(***)$$
where $m$ is the minimal positive integer such that $n\mid mt.$

\vskip10pt

\subsection{} Again using the natural basis of $S_i^l$  one has the following commutative
diagrams in $\A\underline {\mbox{-mod}}$:
\[\xymatrix{\mathcal N(S_{i}^l)\ar[rr]^{\mathcal N(\sigma_{i}^l)} \ar[d]_{\theta_i^l}^-{\wr}
&& \mathcal N(S_{i-1}^{l+1})\ar[d]_{\theta_{i-1}^{l+1}}^-{\wr}\\
S_{i+1-t}^l \ar[rr]^{\sigma_{i+1-t}^l} && S_{i-t}^{l+1};} \ \ \ \ \ \
\ \ \ \ \ \ \ \ \ \ \xymatrix{ \mathcal
N(S_{i}^l)\ar[rr]^{\mathcal N(p_{i}^l)} \ar[d]_{\theta_i^l}^-{\wr}
&&
\mathcal N(S_{i}^{l-1})\ar[d]_{\theta_i^{l-1}}^-{\wr}\\
S_{i+1-t}^l \ar[rr]^{p_{i+1-t}^l} && S_{i+1-t}^{l-1}.}\]

\vskip10pt We justify the commutative diagrams above. Note that  for
any finite quiver $Q$ the bimodule structure of $D(kQ)$ is given
by (using the dual basis)
$$p^*a = \begin{cases} b^*, \ & \mbox{if} \ ab = p, \\ 0, \ &
\mbox{otherwise;}\end{cases} \ \ \ \
\ \ \ \ \ \ \ \mbox{and} \ \ \ \ \ \ \ \ \ \ \ ap^* = \begin{cases} b^*, \ & \mbox{if} \ ba = p, \\
0, \ & \mbox{otherwise.}\end{cases}$$ for any paths $p$ and $a$.
Note that $N(S_{i}^l) = D(\A)\otimes_\A S_{i}^l$ is spanned by
$(\gamma_j^{l'})^*\otimes_\A \gamma_{i+l-t}^{t-u}$, where
$j\in\Bbb Z/n\Bbb Z, \ 1\le l'\le t-1, \ 1\le u\le l.$ By
$(\gamma_j^{l'})^*\otimes_\A \gamma_{i+l-t}^{t-u}=
(\gamma_j^{l'})^*\gamma_{i+l-t}^{t-u}\otimes_\A e_{i+l-t}$ we see
that if $(\gamma_j^{l'})^*\otimes_\A \gamma_{i+l-t}^{t-u}\ne 0$
then $j = i+l-l'-u$; and in this case we have
$$(\gamma_{i+l-l'-u}^{l'})^*\otimes_\A \gamma_{i+l-t}^{t-u}= (\gamma_{i+l-l'-u}^{l'})^*
\gamma_{i+l-t}^{t-u}\otimes_\A e_{i+l-t}
=(\gamma_{i+l-(l'+u)}^{l'+u-t})^*\otimes_\A e_{i+l-t}.$$ This makes
sense only if $l'+u\ge t$. So we have a basis of $N(S_{i}^l)$:
$$(\gamma_{i+l-t-v}^v)^*\otimes_\A e_{i+l-t}, \ \ 0\le v\le l-1.$$
Using the natural basis of $S_{i+1-t}^l$ given in 5.1 we have a
$\A$-isomorphism $\theta_i^l: \ \mathcal N(S_{i}^l)\longrightarrow
S_{i+1-t}^l$ for any $i\in\Bbb Z/n\Bbb Z$ and $1\le l\le t-1$:
$$\theta_i^l: \ (\gamma_{i+l-t-v}^v)^*\otimes_\A e_{i+l-t}\mapsto
\gamma_{i+1+l-2t}^{t-(v+1)}, \ \ 0\le v\le l-1.$$ (One
checks that this is indeed a left $\A$-map.) Note that
$N(\sigma_{i}^l): \ \mathcal N(S_{i}^l)\longrightarrow \mathcal
N(S_{i-1}^{l+1})$ is a natural embedding given by
$$(\gamma_{i+l-t-v}^v)^*\otimes_\A
e_{i+l-t}\mapsto (\gamma_{i+l-t-v}^v)^*\otimes_\A e_{i+l-t}, \ \
0\le v\le l-1;$$ and that $N(p_{i}^l): \ \mathcal
N(S_{i}^l)\longrightarrow \mathcal N(S_{i}^{l-1})$ is a
$\A$-epimorphism given by
$$(\gamma_{i+l-t-v}^v)^*\otimes_\A e_{i+l-t}\mapsto
(\gamma_{i+l-t-v}^{v-1})^*\otimes_\A e_{i+l-1-t}, \ \ 0\le v\le
l-1$$ where $\gamma_{i+l-1-t}^{-1}$ is understood to be $0$. Then
one easily checks the following
$$\sigma_{i+1-t}^l\circ \theta_i^l = \theta_{i-1}^{l+1}\circ \mathcal N(\sigma_i^l); \ \ \ \
p_{i+1-t}^l\circ \theta_i^l = \theta_{i}^{l-1}\circ \mathcal
N(p_i^l).$$ This justifies the commutative diagrams.

Since all these $\theta_i^l$ depend only on $i$ and $l$, which
means that they do not depend on whatever the maps $\sigma_i^l$ or
$p_i^l$ are (this is important for the bi-naturality of a
Calabi-Yau category), it follows, without loss of the generality,
that we can specialize these maps to identities. Thus we have
\begin{align} \mathcal N(S^l_i) = S_{i+1-t}^{l}, \ \ \ \mathcal
N(\sigma^l_i) = \sigma_{i+1-t}^{l}, \ \ \ \mathcal N(p^l_i) =
p_{i+1-t}^{l}.\end{align}

\vskip10pt

\subsection{} By $(***)$ we have $o([1]) <
\infty$, it follows, without loss of generality, that we can
assume $d\ge 0$. For convenience, set $d(t): = 1 +
\frac{(d-1)t}{2}\in \frac {1}{2}\Bbb Z$, whatever $d$ is even or
odd; denote by $N=N(d, n, t)$ the minimal positive integer such
that
\begin{align}
\begin{cases} n\mid N d(t), \ &\mbox{if} \ (d-1)t \ \mbox{is
even}; \\ n\mid N(2d(t)), \ &\mbox{if} \ (d-1)t \ \mbox{is odd}.
\end{cases}\end{align}
(When $2d(t)$ is odd, we will write it together in the following.)
\vskip10pt By $(5.1)$ we have
\ $\Omega^{(2m-1)}(S_i^l) = S_{i+l+(m-1)t}^{t-l}$ \ and \ $
\Omega^{2m}(S_i^l) = S_{i+mt}^{l}$,  hence by
$(5.2)$ we have (remember $G: = \Omega^{d+1}\circ \mathcal N$ and $d\ge 0$)
\begin{align*} G(S_i^l) = S_{i+1+(m-1)t}^{l}, \ \mbox{if} \ d =
2m-1\end{align*} and
\begin{align*} G(S_i^l) = S_{i+l+1+(m-1)t}^{t-l}, \ \mbox{if} \ d = 2m.\end{align*}
Thus, by induction we have
\begin{align} G^{m'}(S_i^l) = S_{i+m'(1+(m-1)t)}^{l} =  S_{i+m'd(t)}^{l}, \ \mbox{if} \ d =
2m-1;\end{align} and
\begin{align} G^{2m'}(S_i^l) = S_{i+m'(2+(2m-1)t)}^{l} = S_{i+m'(2d(t))}^{l}, \ \mbox{if} \ d =
2m,
\end{align}
and
\begin{align} G^{2m'+1}(S_i^l) = S_{i+(m'+1)(2d(t))+l-1-mt}^{t-l}, \ \mbox{if} \ d =
2m,
\end{align}
for $m'\ge 0$.

\subsection{} If $d= 2m-1\ge1$, then by $(5.4)$ we see that
$o(G_{S_i^l}) = N$, where $N$ is as given in $(5.3)$, i.e. $N$ is the
minimal positive integer such that $n\mid N(1+(m-1)t)$. It follows
from Corollary 4.3 and $(5.4)$ that we have

\vskip10pt

\begin{lem} Let $d = 2m-1\ge 1$.
Then $M$ is a minimal $d$-th CY $\A$-module
if and only if $M$ is isomorphic to one of the following

\begin{align} S_i^l\oplus S_{i+d(t)}^l\oplus S_{i+2d(t)}^l\oplus \cdots \oplus
S_{i + (N-1)d(t)}^l, \ \ 1\le l\le t-1, \ \ i\in\Bbb Z/n\Bbb
Z.\end{align}

\vskip10pt In particular, all the minimal $d$-th CY modules
have the same number $N = N(d, n, t)$ of
indecomposable direct summands.
\end{lem}

\vskip10pt

\subsection{} If $d =
2m\ge 0$ and $t$ is odd, then by $(5.6)$ we see
$G^{2m'+1}(S_i^l)\ne S_i^l$ for any $i, l$ (since $t-l\ne l$).
Note that in this case $(d-1)t$ is odd. It follows from $(5.5)$ that
$o(G_{S_i^l}) = 2N$, where $N$ is as given in $(5.3)$, i.e. $N$
is the minimal positive integer such that $n\mid N(2d(t))$. It
follows from Corollary 4.3, $(5.5)$ and $(5.6)$ that we have

\vskip10pt

\begin{lem} Let $t\ge 3$ be an odd integer and $d = 2m\ge 0$. Then $M$ is a
minimal $d$-th CY $\A$-module if and only
if $M$ is isomorphic to one of the following
\begin{align} S_i^l\oplus S_{i+l' + 2d(t)}^{t-l} \oplus
S_{i+2d(t)}^l\oplus S_{i+l'+ 4d(t)}^{t-l} \oplus \cdots \oplus
S_{i + 2d(t)(N-1)}^l\oplus S_{i + l'+ 2d(t)N}^{t-l}, \end{align}
where $l': = l-1-mt,$ \ $1\le l\le t-1$ \ and \ $i\in\Bbb Z/n\Bbb
Z$.

In particular, any minimal $d$-th CY modules has $2N = 2N(d, n, t)$ indecomposable direct
summands.\end{lem}

\vskip10pt

\subsection{} Let $d = 2m\ge 0$ and $t=2s$. Then $d(t) = 1+ (2m-1)s\in\Bbb Z$.

First, we consider $o(G_{S^s_i})$.
In this case $(5.5)$ and $(5.6)$ can be written in a unified way:
\begin{align} G^{m'}(S_i^s) = S_{i+m'd(t)}^{s}, \ m'\ge 0, \ \mbox{if} \ d =
2m, \ t=2s.\end{align} So we have $o(G_{S_i^s}) = N$, where $N$ is
as given in $(5.3)$, i.e. $N$ is the minimal positive integer such
that $n\mid N(1+(2m-1)s)$.

\vskip10pt

Now, we consider $o(G_{S^l_i})$ with $l\ne s, \ 1\le l\le t-1$.
In this case $(5.5)$ and $(5.6)$ are written respectively as:
\begin{align} G^{2m'}(S_i^l) = S_{i+2m'd(t)}^{l}, \ \mbox{if} \ d =
2m, \ t=2s, \ l\ne s,\end{align}
and
\begin{align} G^{2m'+1}(S_i^l) = S_{i+(2m'+1)d(t)+l-s}^{t-l}, \ \mbox{if} \ d =
2m, \ t=2s, \ l\ne s\end{align} for $m'\ge 0.$ Since $l\ne t-l$
for $l\ne s$, it follows that $G^{2m'+1}(S_i^l) \ne S_i^l$. So by
$(5.10)$ we see $o(G_{S_i^l}) = 2N',$ where $N'$ is the minimal
positive integer such that $n\mid 2N'd(t)$. In order to determine
$N'$, we divided into two cases.

{\em Case} 1. \ If $N = N(d, n, t) = N(2m, n, 2s)$ is even, then
$o(G_{S_i^l}) = N$ for any $1\le l\le t-1$. It follows from
Corollary 4.3, $(5.9)$, $(5.10)$, and $(5.11)$ that we have the
following (note that in this case $(5.9)$ is exactly $(5.10)$
together with $(5.11)$, by taking $l=s$)

\vskip10pt

\begin{lem} Let $t = 2s$ and $d = 2m\ge 0$.
Assume that $N = N(d, n, t)$ is even. Then $M$ is a minimal
$d$-th CY $\A$-module if and only if
$M$ is isomorphic to one of the following
\begin{align} S_i^l\oplus S_{i+d(t)+l-s}^{t-l} \oplus
S_{i+2d(t)}^l\oplus S_{i+3d(t)+l-s}^{t-l}\oplus \cdots \oplus S_{i
+ (N-2)d(t)}^l\oplus S_{i + (N-1)d(t)+l-s}^{t-l}\end{align} where
$1\le l\le t-1, \ i\in\Bbb Z/n\Bbb Z$.

\vskip10pt In particular, all the minimal $d$-th CY modules
have the same number $N$ of indecomposable
direct summands.
\end{lem}

\vskip10pt

It remains to deal with

{\em Case} 2. Let $N = N(d, n, t) = N(2m, n, 2s)$ be odd. Since by
definition $N$ is the minimal positive integer such that $n\mid
Nd(t)$, it follows that $N < o(G_{S_i^l}) = 2N' \le 2N.$ It is
easy to see $N'=N$: otherwise $1\le 2N' - N\le N-1$ and $n\mid
(2N' - N)d(t)$, which contradicts the minimality of $N$. It
follows from Corollary 4.3, $(5.9)$, $(5.10)$, and $(5.11)$ that
we have

\vskip10pt

\begin{lem} Let $t = 2s$ and $d = 2m\ge 0$. Assume that $N=N(d, n, t)$ is odd.
Then $M$ is a minimal $d$-th CY $\A$-module if and only if $M$ is isomorphic to one of the following

\begin{align} S_i^s\oplus S_{i+d(t)}^{s} \oplus
S_{i+2d(t)}^s\oplus  \cdots \oplus S_{i +
(N-1)d(t)}^{s}\end{align} where $i\in\Bbb Z_n$, and
\begin{align}\begin{matrix} S_i^l \oplus & S_{i+d(t)+l-s}^{t-l}
\oplus & S_{i+2d(t)}^l  \oplus
& S_{i+3d(t)+l-s}^{t-l} \oplus \cdots  \oplus & S_{i + (N-1)d(t)}^{l}\\
\oplus S_{i+l-s}^{t-l} \oplus & S_{i+d(t)}^{l} \oplus &
S_{i+2d(t)+l-s}^{t-l} \oplus & S_{i+3d(t)}^{l} \oplus  \cdots
\oplus & S_{i + (N-1)d(t)+l-s}^{t-l}\end{matrix}\end{align}
where $l\ne s$, \ $1\le l\le t-1$ and $i\in\Bbb Z_n$.

\vskip10pt In particular, all the  minimal $d$-th CY modules have
either $N$, or $2N$ indecomposable direct summands. \end{lem}

\vskip10pt

\subsection{} By Lemmas 5.1-5.4 all the
minimal $d$-th CY modules of self-injective Nakayama algebras have been classified,
where $d$ is any given integer.
The main result of
this section is as follows.

\vskip10pt

\begin{thm} For any $n\ge 1, \ t\ge 2, \ d\ge 0$, let $N=N(d, n, t)$ be as in $(5.3)$.
Then $M$ is a minimal $d$-th CY \ $\A$-module if and only if $M$
is isomorphic to one of the following

\vskip10pt

$(i)$ \ \ The modules in $(5.7)$, when $d=2m-1$;

$(ii)$ \ \ The modules in $(5.8)$, when $d=2m$ and $t$ is odd;

$(iii)$ \ \ The modules in $(5.12)$, when $d=2m$, $t=2s$, and
$N(d, n, t)$ is even;

$(iv)$ \ \ The modules in $(5.13)-(5.14)$, when $d=2m$, $t=2s$,
and $N(d, n, t)$ is odd.

\vskip10pt

In particular, any minimal $d$-th CY $\A$-module has either $N$ or
$2N$ indecomposable direct summands; and
\begin{align}\operatorname{min} \{d\ge 0 \ |\ N = c(M),  \mbox{or}
\ 2N = c(M) \} \ \le \operatorname{CYdim}(M) < o([1]_M) \le 2n
\end{align} for any minimal Calabi-Yau $\A$-modules $M$, where $c(M)$ is
the number of indecomposable direct summands of $M$.
\end{thm}
\noindent {\bf Proof.}\ \ By Lemmas 5.1-5.4 the CY dimension $d$
of any minimal Calabi-Yau module $M$ satisfies $N(d, n, t) = c(M)
\ \mbox{or} \ \ 2N(d, n, t) = c(M).$ \hfill $\blacksquare$

\vskip10pt

\begin{rem}  The modules in $(5.7)-(5.8)$ and
$(5.12)-(5.14)$ have overlaps. This is because $d$ is not
uniquely determined by a minimal $d$-th CY module. A general
formula of the CY dimensions of the minimal Calabi-Yau modules
seems to be difficult to obtain.

Note that the inequality on the left hand side in $(5.15)$ can not
be an equality in general. For example, take $n = 2, \ t=4, \ m =
2, \ d=2m-1=3$. Then $d(t) = 5$, \  $N = N(3, 2, 4) = 2$, and
$S_i^l\oplus S_{i+1}^l, \ 1\le l\le 3$, are all the minimal $3$-th
CY modules of $CY$ dimension $0$ if $l = 2$, and $1$ if $l=1, \
3$. However, the left hand side in $(5.15)$ is $0$ since $N(0, 2,
4) = 2$. \end{rem}

\vskip10pt

\section{\bf Self-injective Nakayama algebras with indecomposable Calabi-Yau modules}

In this section we determine all the self-injective Nakayama
algebras $\A = \A(n, t), \ n\ge 1, \ t\ge 2$, which admit
indecomposable Calabi-Yau modules.

\vskip10pt

Note that Erdmann and Skowro\'nski have proved in $\S 2$ of [ES]
that self-injective algebras $A$ such that $A\underline
{\mbox{-mod}}$ is Calabi-Yau of CY dimension $0$ (resp. $1$) are
the algebras Morita equivalent to $\A(n, 2)$ for some $n\ge 1$
(resp. $\A(1, t)$ for some $t\ge 3$). So we assume that $t\ge 3$.

\vskip10pt

\begin{thm} \ \ Let $t\ge 3$. Then $\A$ has
an indecomposable Calabi-Yau module if and only if $n$ and $t$
satisfy one of the following conditions

$(i)$ \ $g.c.d. \ (n, \ t) = 1$. This is exactly the case where
$\A\underline {\mbox{-mod}}$ is a Calabi-Yau category. In this
case we have $\operatorname{CYdim}(\A) = 2m-1$, where $m$ is the
minimal positive integer such that $n\mid (m-1)t + 1$.

\vskip10pt

$(ii)$ \ $g.c.d. \ (n,\ t)\ne 1$, \ $t = 2s$, \ and \
$g.c.d. \ (n,\ s) = 1$.  This is exactly the case where $\A\underline
{\mbox{-mod}}$ is not a Calabi-Yau category but admits
indecomposable Calabi-Yau modules.

In this case, we have $g.c.d.\ (n,\ t) = 2$ and

$(a)$\ \  $S_i^s, \ i\in \Bbb Z/n\Bbb Z,$ \ are all the
indecomposable Calabi-Yau modules;

$(b)$ \ \ $S_i^l\oplus S_{i+l-s}^{t-l}, \ 1\le l\le s-1, \
i\in\Bbb Z/n\Bbb Z,$ \ are all the decomposable minimal $2m$-th CY
modules, where $m$ is the minimal non-negative integer such that
$n\mid (2m-1)s + 1$;

$(c)$ \ All of these modules in $(a)$ and $(b)$ have the same CY
dimension $2m$.
\end{thm}

\vskip10pt

{\bf Remark.} Bialkowski and Skowro\'nski [BS] have classified
representation-finite self-injective algebras whose stable
categories are Calabi-Yau. This includes the assertion $(i)$ of
Theorem 6.1.

\vskip10pt

\noindent {\bf Proof.}\ \ If $\A$ has an indecomposable $d$-th CY
module $S_i^l$, then $\mathcal N(S^l_i)
\cong \Omega^{-(d+1)}(S^l_i)$. By
$(5.2)$ and $(5.1)$ we have
\begin{align}1-t\equiv -mt \ (\mbox{mod} \ n), \  \mbox{if} \ d+1 =2m,\end{align}
or
\begin{align} t= t-l, \ 1-t\equiv l-mt \ (\mbox{mod} \ n), \  \mbox{if} \ d+1 =2m-1.\end{align}
In the first case  we get $g.c.d.\ (n, \ t) = 1.$ In the second case we get
$t=2s, \ l=s$ and $g.c.d.\ (n, \ s) = 1.$ Excluding the overlap
situations we conclude that either $g.c.d.\ (n, \ t) = 1;$ \ or \
$g.c.d.\ (n, \ t) \ne 1,  \ t= 2s, \ g.c.d.\ (n, \ s) =
1$.

\vskip10pt

Assume that $g.c.d. \ (n, \ t) = 1$. Then there exists an integer
$m$ such that $n\mid (m-1)t+1$. We chose a positive (otherwise,
add $n(-m+2)t$), and minimal $m$. Set $d: = 2m-1$. Then the same
computation shows that every indecomposable is a $d$-th CY module.
We claim that $\A\underline {\mbox{-mod}}$ is a Calabi-Yau
category. For this, it remains to show $\mathcal N(f) =
\Omega^{-(d+1)}(f) =\Omega^{-2m}(f)$ for any morphism $f$ between
indecomposables of $\A\underline {\mbox{-mod}}$. Since $n\mid
(m-1)t +1$, it follows from $(**)$ in $\S 5$ that
\begin{align}\mathcal N(\sigma_i^l) =
\Omega^{-2m}(\sigma^l_i), \ \   \mathcal N(p_i^l) =
\Omega^{-2m}(p^l_i).\end{align} Since $\A$ is
representation-finite, it follows that $f$ is a $k$-combination of
compositions of irreducible maps $\sigma^l_i$'s and $p^l_i$'s,
hence $\mathcal N(f) = \Omega^{-2m}(f)$ by $(6.3)$. We stress here
that, this argument relies on the fact that all the isomorphisms
$\theta_i^l$ in 5.2 depend only on $i$ and $l$, which means that
they do not depend on whatever the maps $\sigma_i^l$ or $p_i^l$
are. Otherwise we can not take them as identities, and then we can
not get the naturality for the Calabi-Yau category of this case.

This proves the claim, hence $\A\underline {\mbox{-mod}}$ is a
Calabi-Yau category of CY dimension $D$, with $0\le D \le d =
2m-1$. We claim $D=d$. In fact, since every indecomposable is a
$D$-th CY module, and since we have $l$ such that $t\ne t-l$, it
follows from $(6.1)$ and $(6.2)$ that $D = 2m'-1$ with $n\mid
(m'-1)t+1$, hence $D = d$, by the minimality of $m$.

The argument above also proves
that if $\A\underline {\mbox{-mod}}$ is a Calabi-Yau category then
$g.c.d. \ (n, \ t) = 1$.

\vskip10pt

Assume that $g.c.d. \ (n,\ t)\ne 1$, \ $t = 2s$, \ and
$g.c.d. \ (n, \ s) = 1$. Then there exists an integer $m'$ such
that $n\mid (m'-1)s+1$. We choose a positive $m'$. Since $g.c.d.\
(n, \ t) \ne 1$, it follows that $m'$ is even, say $m' =
2m$ with $m\ge 1$. Let $m$ be the the minimal positive integer
such that $n\mid (2m-1)s+1$, and set $d: = 2m$. Then
the same computation shows that \ $S_i^s, \ i\in \Bbb Z/n\Bbb
Z,$ \ are all the indecomposable Calabi-Yau modules.
This proves $(a)$.

By applying Lemma 5.4 to $d$ given above (note that the
corresponding $N = N(n, d, t) = 1$ in this case), we know that $S_i^l\oplus
S_{i+l-s}^{t-l}, \ 1\le l\le t-1, \ l\ne s,  \ i\in\Bbb Z/n\Bbb
Z,$ \ are all the decomposable minimal $2m$-th CY modules. By
symmetry one can consider $1\le l\le s-1:$ if $l>s$ then one can replace
$l$ by $t-l$, since $i = (i+l-s) + (t-l)-s.$ This proves $(b)$.

It remains to prove $(c)$. Let $\operatorname{CYdim}(S^s_i) = d'$.
Since $g.c.d. \ (n,\ t)\ne 1$, it follows
that $d'$ has to be an even integer $2m'\ge 0$ with $n\mid
(2m'-1)s +1$. it follows that $d' = d$, by the minimality of $m$.
Let $\operatorname{CYdim}(S_i^l\oplus S_{i+l-s}^{t-l}) = d'$. Then
we have $\mathcal N(S^l_i) = \Omega^{-(d'+1)}(S_{i+l-s}^{t-l})$.
Since $l\ne s$ it follows from $(5.2)$ and $(5.1)$ that $d'$ has
to be $2m'\ge 0$ with $n\mid (2m'-1)s+1$. Again
by the minimality of $m$ we have $d' = d$. This completes the
proof. \hfill $\blacksquare$

\vskip10pt

\begin{rem} $(i)$ $\operatorname{CYdim}(X)$ usually differs from
$\operatorname{CYdim}(\mathcal A)$ in a Calabi-Yau category
$\mathcal A$.

For example, if $n=3, \ t=4$, then $o([1]) = 6$,
$\operatorname{CYdim}(\A) = 5$, while $\operatorname{CYdim}(S_i^2)
= 2$ and $o([1]_{S_i^2}) = 3, \ \forall \ i\in\Bbb Z/3\Bbb Z$.

However, for $t\ge 3$ and \ $g.c.d.\ (n,\ t) = 1$, if $X$ is
indecomposable and $\operatorname{CYdim}(X)$ is odd then
$\operatorname{CYdim}(X) = \operatorname{CYdim}(\A)$. In fact,
since $o([1]) < \infty$ it follows that $\operatorname{CYdim}(X) =
2m'-1\ge 1$, and $n\mid 1+(m'-1)t$. By Theorem 6.1$(i)$
$\operatorname{CYdim}(\A) = 2m-1$, where $m$ is the minimal
positive integer such that $n\mid 1+(m-1)t$. It follows that
$m'\ge m$, hence $\operatorname{CYdim}(X) \ge
\operatorname{CYdim}(\A)$. On the other hand we have
$\operatorname{CYdim}(X) \le \operatorname{CYdim}(\A)$ by
definition.

$(ii)$ \ Consider the algebra $A(t): = kA_\infty ^\infty/J^t$ where
$A_\infty^\infty$ is the infinite quiver
$$ \cdots \longrightarrow \bullet
\longrightarrow \bullet \longrightarrow \bullet \longrightarrow
\bullet \longrightarrow \bullet\longrightarrow \cdots.$$ Then
$A(t)\underline {\mbox{-mod}}$ has a Serre functor, and there is a
natural covering functor $A(t)\underline {\mbox{-mod}}$
$\longrightarrow \A(n, t)\underline {\mbox{-mod}}$ ([Gab], 2.8).
But one can prove that in any case $A(t)\underline {\mbox{-mod}}$
is not a Calabi-Yau category.
\end{rem}

\vskip10pt

{\bf Acknowledgements.} This work is done during a visit of the
second named author at Universit\'e de Montpellier 2, supported by
the CNRS of France. He thanks the first named author and his group
for the warm hospitality, the D\'epartement de Math\'ematiques of
Universit\'e de Montpellier 2 for the working facilities, and the
CNRS for the support. We thank Bernhard Keller for helpful
conversations.

\vskip30pt

\bibliography{}

\end{document}